\documentclass[12pt,reqno]{amsart}
\usepackage{amsmath,amsthm,amsfonts,amssymb,amscd}
\usepackage[latin1]{inputenc}
\usepackage{psfrag,pslatex}
\usepackage{epsfig}
\usepackage{a4wide}

\numberwithin{equation}{section}

\newcommand{\ov}{\overline}
\newcommand{\un}{\underline}

\newcommand{\dist}{\operatorname{dist}}

\newcommand{\inter}{\operatorname{int}}

\newcommand{\diam}{\operatorname{diam}}
\newcommand{\supp}{\operatorname{supp}}
       
\newcommand{\be} {\beta}        
       
\newcommand{\de} {\delta}       
\newcommand{\ep} {\varepsilon}

\renewcommand{\th} {\theta}

\newcommand{\vfi}{\varphi}

\def \RR {{\mathbb R}}
\def \BB {{\mathcal B}}

\def \AA {{\mathcal A}}
\def \CC {{\mathcal C}}

\def \NN {{\mathbb N}}
\def \SS {{\mathbb S}}

\def \EE {{\mathbb E}}

\def \TT {{\mathbb T}}
\def \cT {{\mathcal T}}

\newcommand{\dem}{\begin{proof}}
\newcommand{\cqd}{\end{proof}}

\newcommand{\qand}{\quad\text{and}\quad}

\newcommand{\cP}{{\mathcal P}}

\newcommand{\cB}{{\mathcal B}}

\newtheorem{maintheorem}{Theorem}

\newcommand{\cmt}{\begin{maintheorem}}
\newcommand{\fmt}{\end{maintheorem}}

\newtheorem{maincorollary}[maintheorem]{Corollary}

\newcommand{\cmc}{\begin{maincorollary}}
\newcommand{\fmc}{\end{maincorollary}}

\newtheorem{T}{Theorem}[section]
\newcommand{\cte}{\begin{T}}
\newcommand{\fte}{\end{T}}

\newtheorem{Corollary}[T]{Corollary}
\newcommand{\cco}{\begin{Corollary}}
\newcommand{\fco}{\end{Corollary}}

\newtheorem{Proposition}[T]{Proposition}
\newcommand{\cpr}{\begin{Proposition}}
\newcommand{\fpr}{\end{Proposition}}

\newtheorem{Lemma}[T]{Lemma}
\newcommand{\cle}{\begin{Lemma}}
\newcommand{\fle}{\end{Lemma}}

\newtheorem{Definition}{Definition}
\newcommand{\cde}{\begin{Definition}}
\newcommand{\fde}{\end{Definition}}

\newcommand{\csle}{\begin{Lemma}}
\newcommand{\fsle}{\end{Lemma}}

\theoremstyle{remark}

\newtheorem{Remark}[T]{Remark}
\newcommand{\cre}{\begin{Remark}}
\newcommand{\fre}{\end{Remark}}

\begin{document}

\author{Vitor Araujo}

\address{Centro de Matemática da
  Universidade do Porto, Rua do Campo Alegre 687, 4169-007
  Porto, Portugal} \email{vdaraujo@fc.up.pt}
\urladdr{http://www.fc.up.pt/cmup/vdaraujo}

\thanks{Trabalho parcialmente financiado por CMUP. O CMUP é
  financiado por FCT, no âmbito de POCTI-POSI do Quadro
  Comunitário de Apoio III (2000-2006), com fundos
  comunitários (FEDER) e nacionais.  A.T. would like to
  thank FCT for financial support and CMUP for warm
  hospitality. V.A. was also partially supported by grant
  FCT/SAPIENS/36581/99.}

\author{Ali Tahzibi} 

\address{Departamento de Matem\'atica,
  ICMC-USP São Carlos, Caixa Postal 668, 13560-970 São
  Carlos-SP, Brazil.} 
\email{tahzibi@icmc.sc.usp.br}
\urladdr{http://www.icmc.sc.usp.br/~tahzibi}

\dedicatory{Dedicated to C. Gutierrez on the occasion of his
  60th birthday.}

\keywords{Intermittent maps, stochastic stability,
   equilibrium states, random perturbations, physical measures}

\subjclass{Primary: 37D25. Secondary: 37D30, 37D20.}

\renewcommand{\subjclassname}{\textup{2000} Mathematics Subject Classification}

\date{\today}

\setcounter{tocdepth}{1}

\title{Stochastic stability at the boundary of expanding
  maps}

\begin{abstract}
  We consider endomorphisms of a compact manifold which are
  expanding except for a finite number of points and prove
  the existence and uniqueness of a physical measure and
  its stochastical stability. We also characterize the
  zero-noise limit measures for a model of the intermittent
  map and obtain stochastic stability for some values of the
  parameter. The physical measures are obtained as
  zero-noise limits which are shown to satisfy Pesin´s
  Entropy Formula.
\end{abstract}

\maketitle

\tableofcontents

\section{Introduction}

After the long and deep developments in the last decades on
the structural stability theory of dynamical systems, we
know that this form of stability is too strong to be a
generic property. Recently there has been some emphasis on
the study of stochastic stability of dynamical systems,
among other forms of stability.

 On the one hand, one of the challenging problems of smooth
 Ergodic Theory is to prove the existence of "nice"
 invariant measures called physical measures or sometimes
 SRB (Sinai-Ruelle-Bowen) measures. On the other hand, a
 natural formulation of stochastic stability of dynamical
 systems assumes the existence of physical measures.
 However, the characterization of zero-noise limit measures
 involved in the study of stochastic stability may provide
 ways to construct physical measures.  In this work the
 study of zero-noise limit measures for endomorphisms which
 are expanding except at a finite number of points yields a
 construction of physical measures and also their stochastic
 stability.

 Let $M$ be a compact and connected Riemannian manifold and
 $\cT:=C^{1+\alpha}(M,M)$ be the space of $C^{1+\alpha}$
 maps of $M$ where $\alpha>0$. We write $m$ for some fixed
 measure induced by a normalized volume form on $M$ that we
 call \emph{Lebesgue measure}, $\dist$ for the Riemannian
 distance on $M$ and $\|\cdot\|$ for the induced Riemannian
 norm on $TM$.
 
 We recall that an invariant probability measure $\mu$ for a
 transformation $T:M\to M$ on a manifold $M$ is
 \emph{physical} if the \emph{ergodic basin}
\[
B(\mu)=\left\{x\in M:
\frac1n\sum_{j=0}^{n-1}\varphi(T^j(x))\to\int\varphi\,
d\mu\mbox{  for all continuous  } \varphi: M\to\RR\right\}
\]
has positive Lebesgue measure.

Let $(\theta_\ep)_{\ep>0}$ be a family of Borel probability
measures on $(\cT, \BB(\cT))$, where we write $\BB(X)$ the
Borel $\sigma-$algebra of a topological space $X$. We are
dealing with random dynamical systems generated by
independent and identically distributed maps of $\cT$ and
$\theta_\ep$ will be the common probability distribution when
choosing the maps to generate random dynamics.

We say that a probability measure $\mu^\ep$ on $M$ is
\emph{stationary for the random system $(\hat T,
  \theta_\ep)$} if the following holds
\begin{equation}
  \label{eq:1}
  \int\!\int \varphi(T(x)) \,d\mu^\ep(x)d\theta_{\ep}(T) =
  \int \varphi \,d\mu^\ep \quad \mbox{for all
  continuous  } \varphi:M\to\RR.
\end{equation}

We assume that the support of $\theta_\ep$ shrinks to $T$
when $\ep \to 0$ in a suitable topology.  A classical result
in random dynamical systems (see~\cite{Ki88} or~\cite{Ar00})
implies that \emph{every weak$^*$ accumulation point of the
  stationary measures $(\mu^\ep)_{\ep>0}$ when $\ep \to 0$
  is a $T$-invariant probability measure,} which is called a
\emph{zero-noise limit measure.}  This naturally leads to
the study of the kind of zero noise limits that can arise
and to the notion of stochastic stability.

\begin{Definition}
  \label{def.stocasticstable}
  A map $T$ is stochastically stable (under the random
  perturbation $(\hat T,\theta_\ep)_{\ep>0}$) if every
  accumulation point $\mu$ of the family of stationary
  measures $(\mu^\ep)_{\ep>0}$, when $\ep\to0$, is a linear
  convex combination of the physical measures of $T$.
\end{Definition}

Uniformly expanding maps and uniformly hyperbolic systems
are known to be stochastically
stable~\cite{Ki86a,Ki88,Vi97b,Yo85}. Some non-uniformly
hyperbolic systems, like quadratic maps, Hénon maps and
Viana maps, were shown to be stochastically stable much more
recently~\cite{AA03,BaV96,BeV2}. These systems either
exhibit expansion/contraction everywhere or are
expanding/contracting away from a critical region with slow
recurrence rate to it. This allows for a probabilistic
argument which shows that the visits to a neighborhood of
the critical region are negligible on the average, and also
that this behavior persists under small random
perturbations.

It is not obvious how to apply the standard techniques to
systems whose typical orbits do not have a slow recurrence
rate of visits to the non-hyperbolic regions. This is the
case of \emph{intermittent maps}~\cite{manneville1980}.
These applications are expanding, except at a neutral fixed
point. The local behavior near this neutral point is
responsible for various phenomena.

Consider $\alpha>0$ and the map $T:
[0,1] \rightarrow [0,1]$ defined as follows

\begin{equation} \label{intermit}
T(x)= \left\{
\begin{array}{ll}
 x + 2^{\alpha} x^{1+\alpha} \qquad &  x \in [0,\frac{1}{2})\\
 x - 2^{\alpha} (1-x)^{1 + \alpha}  &  x \in [\frac{1}{2},1]
\end{array} \right.
\end{equation}
This map defines a $C^{1+\alpha}$ map of the unit circle
$\SS$ into itself.  The unique fixed point is $0$ and
$DT(0)= 1.$ The above family of maps provides many
interesting results in Ergodic Theory. If $\alpha \geq 1$,
i.e if the order of tangency at zero is high enough, then
the Dirac mass at zero $\delta_0$ is the unique physical
probability measure and so the Lyapunov exponent of Lebesgue
almost all points vanishes~\cite{thaler1983}.  The situation
is completely different for $0 < \alpha < 1$: in this case
there exists a unique absolutely continuous invariant
probability measure $\mu_{SRB}$, which is therefore a
physical measure and whose basin has full Lebesgue
measure~\cite{thaler1980}.

Another point of interest is that these maps provide
examples of dynamical systems with polynomial decay of
correlations. M.  Holland has obtained even sub-polynomial
rate of mixing modifying these intermittent maps
\cite{Ho03}. In particular $\mu_{SRB}$ is always mixing,
when it exists.

\subsection{Statement of the results}
\label{sec:statement-results}

We consider additive noise applied to a map $T$ of $\SS$ with an
indifferent fixed point at $0$ and expanding everywhere else, as in
example \eqref{intermit}.  Let $ \alpha >0$ be fixed and consider $T_t
:= T + t$ for $|t| \leq \ep.$ Then $\hat T:[-1/2,1/2]\to
C^{1+\alpha}(\SS,\SS), t\mapsto T_t$ is a (smooth) family of
$C^{1+\alpha}$ maps of $\SS$.

Let $\theta_\ep$ be an absolutely continuous probability measure, with
respect to the Lebesgue measure $m$ on $\SS$, whose support is
contained in $[-\ep, \ep]$ (e.g. $\theta_{\ep} = (2\ep)^{-1} m \mid
[-\ep,\ep], \ep>0$).  This naturally induces a probability measure on
$\cT=\{T_t,t\in[-1/2,1/2]\}$ which we denote by the same symbol
$\th_\ep$ (the meaning being clear from the context).

In this setting it is well known that there always exist a
stationary probability measure $\mu^\ep$ for all $\ep>0$.
Moreover this measure is ergodic and is the unique
absolutely continuous stationary measure for $(\hat T,
\theta_\ep)$ (see Subsection \ref{sec:topological-mixing}).

Let us fix now $\alpha\in(0,1)$ and let
\[
\EE = \{ t \delta_0 + (1-t) \mu_{SRB} : 0 \leq t \leq 1 \}
\]
be the set of linear convex combinations of the Dirac mass
at $0$ with the unique absolutely continuous invariant
probability measure for these maps.

\cmt
 \label{th.equil}
 Let $\mu_0$ be any accumulation point of the stationary
 measures $(\mu^{\ep})_{\ep>0}$ when $\ep\to0$ for
 the random perturbation $(\hat
 T,\theta_\ep)_{\ep>0}$ with $\alpha\in(0,1)$.
 Then $\mu_0 \in \EE$.  \fmt
 
 In the case $\alpha\ge1$ there does not exist any
 absolutely continuous invariant probability measure.
 However the Dirac measure $\delta_0$ is the unique physical
 measure. In this case we are able to obtain stochastic
 stability.

\cmt
\label{th.stocstable}
Let $\alpha\ge1$
in~\eqref{intermit} and let $(\mu^\ep)_{\ep>0}$ be
the family of stationary measures for the random
perturbation $(\hat T,\theta_\ep)_{\ep>0}$. Then
$\mu^\ep\to\delta_0$ when $\ep\to0$ in the
weak$^*$ topology.
\fmt

\begin{figure}[htbp]
  \centering
  \includegraphics[width=12cm, height=6cm]{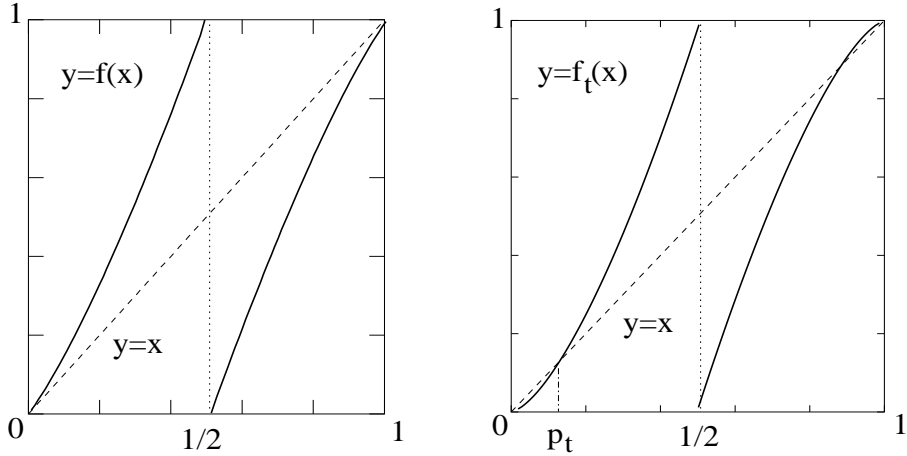}
  \caption{The map $f=T$ (left) and the map $f_t$ for $0<t<1$ (right).}
  \label{fig1}
\end{figure}

However, taking a different family $f_t$ unfolding the
saddle-node at $0$, e.g.
\begin{equation}
  \label{eq:nonstochfamily}
  f_t(x)=\left\{
\begin{array}{ll}
tx + 2^{\alpha}(2-t) x^{1+\alpha} \qquad &  x \in [0,\frac{1}{2})\\
1-t(1-x) - 2^{\alpha}(2-t) (1-x)^{1 + \alpha}  
&  x \in [\frac{1}{2},1]
\end{array} \right.
\end{equation} 
with $\alpha\in(0,1)$ we obtain an example of \emph{non-stochastic
stability}. In fact, since $f'_t(0)=t$ then for $t<1$ the fixed point
$0$ is a sink for $f_t$ (see Figure~\ref{fig1}) and we prove that the
physical measure for the random system is always $\de_0$ for
restricted choices of the probability measures $\th_\ep$.

\cmt\label{thm.non-stochstable} For every small enough $\ep>0$ there
are $a(\ep)<b(\ep)<1$ such that $a(\ep)\to1$ when $\ep\to0$ and, for
any given probability measure $\th_\ep$ supported in
$[a(\ep),b(\ep)]$, the unique stationary measure $\mu^\ep$ for the
random system equals $\de_0$.  \fmt

Since $f_1=T$ admits an absolutely continuous invariant
measure $\mu_{SRB}$ and clearly $\de_0$ cannot converge to
this physical measure, we have an example of a
stochastically unstable system (under this kind of
perturbations).

Using the same kind of additive perturbations considered in
Theorems~\ref{th.equil} and~\ref{th.stocstable},
our methods provide the following results for maps in higher
dimensions.

\cmt
\label{th.limitSRB}
Let $f:M\to M$ be a $C^{1+\alpha}$ local diffeomorphism such that
\begin{enumerate}
\item $\|Df(x)^{-1}\|\le1$ for all $x\in M$;
\item $K=\{x\in M: \|Df(x)^{-1}\|=1\}$ is finite and $|\det
  Df(x)|>1$ for every $x\in K$.
\end{enumerate}
Then, for any non-degenerate random perturbation $(\hat
f,\th_\ep)_{\ep>0}$, there exists a unique
ergodic stationary probability measure $\mu^\ep$ for all
$\ep>0$.  Moreover $\mu^\ep$ converges, in the weak$^*$
topology when $\ep\to0$, to a unique absolutely continuous
$f$-invariant probability measure $\mu_0$ whose basin has
full Lebesgue measure, and $f$ is stochastically stable.
\fmt

Here we will assume that $M$ is a $n$-dimensional torus
since the maps $f$ satisfying the conditions on
Theorem~\ref{th.limitSRB} are at the boundary of expanding
maps, which can only exist on special
manifolds~\cite{Sh69,Gr81}, the best known example being the
tori. Since these manifolds are parallelizable, we can
define additive perturbations just as we did on the circle.
If $\TT^n$ is a $n$-dimensional torus, then
$T\TT^n\simeq\RR^n$ and $\hat f: B\subset \RR^n\to
C^{1+\alpha}(M,M), v\to f+v$, where $B$ is a ball around the
origin of $\RR^n$ (together with a family
$(\th_\ep)_{\ep>0}$ of absolutely continuous probability
measures on $B$, see Subsection~\ref{sec:non-deg-pert} for
the definition of non-degenerate random perturbation) will
be a the kind of additive perturbation we will consider.

These results will be derived from the following more
technical one, but also interesting in itself.

\cmt\label{th.limitequilibrium}
Let $f:M\to M$ be a $C^{1+\alpha}$ local diffeomorphism such that
\begin{enumerate}
\item $\|Df(x)^{-1}\|\le1$ for all $x\in M$;
\item $K=\{x\in M: \|Df(x)^{-1}\|=1\}$ is finite.
\end{enumerate}
Then, for any non-degenerate random perturbation $(\hat
f,\th_\ep)_{\ep>0}$, every weak$^*$ accumulation
point $\mu$ of the sequence $(\mu^\ep)_{\ep>0}$, when
$\ep\to0$, is an equilibrium state for the potential $-\log |\det
Df(x)|$, i.e.
\begin{equation}
  \label{eq:3}
  h_\mu(f) = \int \log |\det Df(x)| \, d\mu (x).
\end{equation}
Moreover every equilibrium state $\mu$ as above is a convex
linear combination of an absolutely continuous invariant
probability measure with finitely many Dirac measures
concentrated on periodic orbits whose Jacobian equals 1.
\fmt

Cowieson and Young have presented results similar to ours
for $C^2$ or $C^\infty$ diffeomorphisms. However their
assumptions are on the convergence of the sum of the
positive Lyapunov exponents for the random maps to the same
sum for the original map, and they obtain \emph{SRB
  measures}, not necessarily physical ones,
see~\cite{CoYo2004} for more details. We make much stronger
assumptions on both the kind of maps being perturbed
(expanding except at finitely many points) and the kind of
perturbations used (additive besides being absolutely
continuous), and we obtain physical measures for $C^{1+\alpha}$
endomorphisms.

In what follows, we first present some examples of
applications and then general results about random dynamical
systems (Section~\ref{sec:preliminary-results}) to be used
to prove Theorem~\ref{th.limitequilibrium}
(Section~\ref{sec:zero-noise-limits}).  At this point we are
ready to obtain Theorem~\ref{th.limitSRB}
(Section~\ref{sec:exist-acim-stoch}).  Finally we apply the
ideas to the specific case of the intermittent maps (Section
\ref{sec:rand-pert-interm}), completing the proof of
Theorems~\ref{th.equil},~\ref{th.stocstable}
and~\ref{thm.non-stochstable}.


\subsection{Examples}
\label{sec:examples-1}

In what follows we write $\TT$ for $\SS\times\SS$ and
consider $\SS=[0,1]/\{0\sim 1\}$. We always assume that
these spaces are endowed with the metrics induced by the
standard Euclidean metric through the identifications. The
Lebesgue measure on these spaces will be denoted by $m$
(area) on $\TT$ and $m_1$ (length) on $\SS$.

An extra example is the intermittent map itself, dealt with
in Section~\ref{sec:rand-pert-interm}.

\subsubsection{Direct product
  ``intermittent$\times$expanding''}
\label{sec:direct-prod-interm}

Let $f:\TT\to\TT, (x,y)\mapsto (T_\alpha(x),g(y))$, where
$T_\alpha$ is defined at the Introduction with $\alpha>0$,
and $g:\SS\to\SS$ is $C^{1+\alpha}$, admits a fixed point
$g(0)=0$ and $g'=Dg>1$.

Since $f$ is a direct product, if $\alpha\in(0,1)$, then $f$
admits an invariant probability measure
$\nu=\mu_{\alpha}\times\lambda$, where $\mu_\alpha$ is the unique absolutely
continuous invariant measure for $T_\alpha$ and $\lambda$ is
the unique absolutely continuous invariant measure for $g$,
i.e., $\mu_\alpha\ll m_1$ and $\lambda\ll m_1$. Hence the product
measure is absolutely continuous: $\nu\ll m=m_1\times m_1$.
These measures are ergodic and also mixing, and the basins of $\mu_\alpha$ and
$\lambda$ equal $\SS, m_1\bmod0$. Thus their direct product
$\nu$ is ergodic and so $B(\nu)=\TT, m\bmod 0$.

If $\alpha\ge1$, then $\nu=\de_0\times\lambda$ is again an
ergodic invariant probability measure for $f$ with
$B(\nu)=\TT, m\bmod 0$, since $\lambda$ is the same as
before and so is mixing for $g$, and $\de_0$ is
$T_\alpha$-ergodic, with the basin of both measures equal to
$\SS$.

Here $K=\{0\}\times\SS$ (the definition of $K$ is given at
the statement of Theorem~\ref{th.limitSRB}) is not finite,
and the conclusion of Theorem~\ref{th.limitSRB} does not
hold when $\alpha\ge1$: we have a physical measure which is
not absolutely continuous with respect to $m$. Note that
clearly $\|(Df)^{-1}\|\le1$ everywhere and since $K$
contains fixed (and periodic) points, $f$ is not uniformly
expanding.

\subsubsection{Direct product
  ``intermittent$\times$intermittent''}
\label{sec:direct-prod-2}

Let $f:\TT\to\TT, (x,y)\mapsto (T_\alpha(x),T_\beta(y))$
where $\alpha,\beta>0$. Now
$K=\{0\}\times\SS\cup\SS\times\{0\}$ and, by the same
reasoning of the previous example, the probability measure
$\nu=\mu_\alpha\times\mu_\beta$ is the unique physical
measure for $f$. Moreover $B(\nu)=\TT$ as before. However $\nu$
is absolutely continuous with respect to $m$ if, and only
if, $\alpha,\beta\in(0,1)$.

\subsubsection{Skew-product
  ``intermittent$\rtimes$expanding''}
\label{sec:skew-prod-interm}

Let $f:\TT\to\TT, (x,y)\mapsto (T_\alpha(x)+\eta y,
g(y))$, for $\alpha>0$, $\eta\in(0,1)$,
and $g:\SS\to\SS$ as in
example~\ref{sec:direct-prod-interm}.

In this case we easily calculate $Df=\left(
\begin{array}{cc}
DT_\alpha  &  \eta
\\
0 & Dg
\end{array}
\right)$
and so $K=\{(0,0)\}$.

Clearly $\|(Df)^{-1}\|\le1$ everywhere and since $K$ is a
fixed point the map $f$ is not uniformly expanding. Applying
Theorem~\ref{th.limitSRB} we get an absolutely continuous
invariant probability measure $\mu$ for $f$ with
$\int\log\|(Df)^{-1}\|\,d\mu<0$. Hence the Lyapunov
exponents for Lebesgue almost every point on the basin of
$\mu$ are all positive, so $f$ is a non-uniformly expanding
transformation.

This map is stochastically stable,
since every weak$^*$ accumulation point of
$(\mu^\ep)_{\ep>0}$ when $\ep\to0$ equals $\mu$ by the
uniqueness part of Theorem~\ref{th.limitSRB}.
We stress that since the value of $\alpha$ played no role in
the arguments, these conclusions hold for any $\alpha>0$.


\section{Preliminary results}
\label{sec:preliminary-results}

Throughout this section we outline some general results
about random dynamical systems to be used in what follows.

Having a parameterized family of maps $\hat T: X\to \cT,
t\mapsto T_t$, where $X$ is some connected compact metric
space, enables us to identify a
sequence $T_1,T_2,\dots$ of maps from $\cT$ with a sequence
$\omega_1,\omega_2,\dots$ of parameters in $X$. The
probability measure $\theta_\ep$ can then be assumed to be
supported on $X$.

We set $\Omega=X^\NN$, the space of sequences
$\omega=(\omega_i)_{i\ge1}$ with elements in $X$. Then we
endow $\Omega$ with the standard infinite product topology,
which makes $\Omega$ a compact metrizable space, with
distance given by (for example)
$d(\omega,\omega')=\sum_{j\ge1} 2^{-1}
  d_X(\omega_j,\omega'_j)$ where $d_X$ is the distance on
  $X$. We also take the standard product probability measure
  $\theta^\ep=\theta_\ep^\NN$, which makes
  $(\Omega,\cB,\theta^\ep)$ a probability space. Here
  $\cB=\BB(X)$ is the $\sigma$-algebra generated by cylinder sets,
  that is, the minimal $\sigma-$algebra of subsets of
  $\Omega$ containing all sets of the form $\{ \omega \in
  \Omega : \omega_1 \in A_1, \omega_2 \in A_2, \cdots ,
  \omega_l \in A_l \}$ for any sequence of Borel subsets
  $A_i \subset X, i= 1, \cdots , l$ and $l\ge1$.

The following skew-product map is the natural setting for
many definitions connecting random with standard dynamical
systems
\begin{eqnarray*}
   S &:& \Omega\times M  \to \Omega\times M \\
   & & (\omega,x) \mapsto (\sigma(\omega),T_{\omega_1}(x))
\end{eqnarray*}
where $\sigma$ is the left shift on $\Omega$, defined as
$(\sigma(\omega))_n = \omega_{n+1}$ for all $n\ge1$. It is
an exercise to check that $\mu^\ep$ is a stationary measure for
the random system $(\hat T,\theta_\ep)$ (i.e. satisfying
\eqref{eq:1}) if, and only if, $\theta^\ep \times \mu^\ep$ on
$\Omega \times M$ is invariant by $S.$
 Ergodicity of
stationary measures is defined in a natural way. A Borel
set $A \in \BB(M)$ is called invariant if for $\mu^\ep$-almost
every point $x \in M$
 \begin{eqnarray*}
  x \in A &\Rightarrow& T_t(x) \in A \quad \text{   for }
  \theta_\ep-\text{almost every   } t \in X; \mbox{ and} \\
  x \in A^c &\Rightarrow& T_t(x) \in A^c \quad \text{for }
  \theta_\ep-\text{almost every   }  t \in X.
 \end{eqnarray*}

\cde \label{ergodicity} A stationary measure $\mu^\ep$ is said
to be ergodic if every Borel invariant set has either
$\mu^\ep$-measure zero or one.  \fde

It is not difficult to prove that $\mu^\ep$ is ergodic if and
only if $\theta^{\ep} \times \mu^\ep$ is an ergodic measure for
$S$ (see for example \cite{LQ95}).

\subsection{Metric entropy of Random Dynamical Systems}
\label{sec:metr-entr-rand}

The notion of metric entropy can be defined for random
dynamical systems in different ways. We point out two
definitions which will be used in this paper and relate
them.  The following results  can be found in the book of
Kifer \cite[Section II]{Ki86}.

Let $\mu$ be a stationary measure for the random system
$(\hat T,\theta_\ep)$ as defined in the beginning of this
section.

\begin{T}{\cite[Thm. 1.3]{Ki86}} \label{thm.metr-entr-rand}
 For any finite measurable partition $\xi$ of $M$ the
  limit
 $$
 h_{\mu^\ep}((\hat T,\theta_\ep), \xi) = \lim_{n \rightarrow \infty}
 \frac{1}{n} \int H_{\mu^\ep} \big(
 \bigvee_{k=0}^{n-1}(T^k_{\omega})^{-1} \xi \big) d
 \theta^{\ep} (\omega)
 $$
 exists. This limit is called the entropy of the random
 dynamical system with respect to $\xi$ and to $\mu^\ep$.
\end{T}

\cre\label{re.infimum}
As in the deterministic case the above limit can be
replaced by the infimum.
\fre

\cde \label{def.entropy} The metric entropy of the random
dynamical system $(\hat T,\theta_\ep)$ is given by
$h_{\mu^\ep} (\hat T, \theta_\ep)= \sup h_{\mu^\ep} ((\hat
T,\theta_\ep), \xi)$, where the supremum is taken over all
measurable partitions.  \fde

It seems natural to define the entropy of a random system by
$h_{\theta^{\ep} \times \mu^\ep} (S)$ where $S$ is the
corresponding skew-product map.  Kifer \cite[Thm. 1.2]{Ki86} shows
that this definition is not very convenient: under some mild
conditions the entropy of $S$ is infinite. However considering an
appropriate $\sigma-$algebra, the conditional entropy of
$\theta^{\ep} \times \mu^\ep$ coincides with the entropy as
defined in Definition \ref{def.entropy}.

Let $\BB \times M$ denote the minimal $\sigma-$algebra
containing all products of the form $A \times M$ with $A \in
\BB.$ In what follows we denote by $h_{ \theta^{\ep}\times
  \mu^\ep}^{\BB \times M} (S)$ the conditional metric
entropy of $S$ with respect to the $\sigma$-algebra $\BB
\times M.$ (See e.g.~\cite{Bi65} for a definition and
properties of conditional entropy.)

\begin{T}{\cite[Thm. 1.4]{Ki86}}
\label{thm.randentropyS}
Let $\mu^\ep$ be a stationary probability measure for the
random system $(\hat T,\theta_\ep)$. Then
$h_{\mu^\ep}(\hat T,\theta_\ep) = h_{
  \theta^{\ep} \times \mu^\ep}^{ \BB \times M} (S)$.
\end{T}

The useful Kolmogorov-Sinai result about generating partitions is
also available in a random version. We denote $\AA=\BB(M)$ the Borel
$\sigma$-algebra of $M$.

\begin{T}{\cite[Cor. 1.2]{Ki86}}
  \label{thm.KSrandom}
  If $\xi$ is a random generating partition for $\AA$, that
  is $\xi$ is a finite partition of $M$ such that
  $$
  \bigvee_{k=0}^{+\infty} (T_{\omega}^k)^{-1} \xi = \AA \quad
  \text{for} \quad \theta^{\ep}-\text{almost all } \omega
  \in \Omega,
  $$
  then $h_{\th^\ep\times\mu^\ep}(\hat
  T,\theta_\ep)=h_{\mu^\ep}((\hat
  T,\theta_\ep),\xi)$.
\end{T}


\subsection{Topological mixing}
\label{sec:topological-mixing}

Here we show that in the setting of
Theorems~\ref{th.limitSRB} and \ref{th.limitequilibrium} we
always have topological mixing for the transformation $T$.
\emph{This ensures uniqueness of stationary measures under
non-degenerate random perturbations}, as we shall see.

Since $T:M\to M$ is a local diffeomorphism on a compact Riemannian
manifold, there is a positive number $\rho$ such that $T\mid
B(x,\rho)$ is a diffeomorphism onto its image and $B(x,\rho)$ is a
convex neighborhood for every $x\in M$, i.e., for every pair of points
$y,z$ in $B(x,\rho)$ there exists a smooth geodesic $\gamma:[0,1]\to
M$ connecting them whose length equals $\dist(y,z)$ and, moreover,
$\gamma\mid [s,t]$ is the curve of minimal length between any pair of
points $\gamma(s),\gamma(t)$ with $s<t,s,t\in [0,1]$.

\cle\label{lem.topmix}
Let $T:M\to M$ satisfy the conditions of
Theorem~\ref{th.limitequilibrium}. Then for every open
subset $U$ there exists an iterate $n\ge1$ such that $T^n(U)=M$.
\fle

\dem

Arguing by contradiction, let us suppose that there exists a
ball $B(x,r)$, for some $x\in M$ and small $r>0$, such that
$T^n(B(x,r))\neq M$ for every $n>1$.

In what follows we fix $n>1$ such that $T^k(B(x,r))\neq M$
for all $k=1,\dots,n$. Then there is $y\in M\setminus
T^n(B(x,r))$ and a smooth curve $\gamma:[0,1]\to M$ such
that $\gamma(0)=T^n(x), \gamma(1)=y$ and whose length is
bounded by $\kappa=\diam (M)+1$ (which is finite, because
$M$ is compact).

Now we fix $\delta\in(0,r/(10k))$ and $\delta'\in(0,\delta)$
small enough such that
\begin{itemize}
\item the $k$ connected components of $B(K,\de')$ are convex
  neighborhoods;
\item the connected components of $T(B(K,\de'))$ (there are
  at most $k$ of them) are also
  convex neighborhoods with diameter smaller than $2\de$.
\end{itemize}
Moreover we choose $\lambda_1\in(0,1)$ such that
$r>2k\delta/\lambda_1$ and set $\lambda=\max\{\|DT(x)^{-1}\| : x\in
M\setminus B(K,\delta')\}$ and
$\lambda_0=\lambda+\lambda_1(1-\lambda)<1$.

We write $\gamma_0$ to denote a smooth curve such that
$T\circ\gamma_0=\gamma$ in what follows.

Let $\ov\gamma_0=\gamma_0\mid \gamma_0^{-1} B(K,\de')$ be the
portion of $\gamma_0$ inside $B(K,\de')$. Since every
connected component of $\ov\gamma=T\circ\ov\gamma_0$ is
inside a convex neighborhood of diameter at most $2\de$, we
may assume that the length $\ell(\ov\gamma)$ of $\ov\gamma$
is at most $2\de$. For otherwise we may replace $\ov\gamma$
by portions of minimizing geodesics connecting the endpoints
of each connected component, with smaller total length.  Now
we obtain, by the non-contracting character of $T$ and by
the definitions of the constants above
\begin{eqnarray*}
  \ell(\gamma_0)
  &=&
  \ell(\gamma_0\setminus\ov\gamma_0)+\ell(\ov\gamma_0)
  \le \lambda\cdot\ell(\gamma\setminus\ov\gamma)+\ell(\ov\gamma)
  \\
  &=&
  \ell(\gamma)\left(
    \frac{\lambda(\ell(\gamma)-\ell(\ov\gamma)) +
  \ell(\ov\gamma)}{\ell(\gamma)}
    \right)
    =
    \ell(\gamma)\left(
      \lambda + (1-\lambda)\frac{\ell(\ov\gamma)}{\ell(\gamma)}
      \right)
  \\
  &\le&
  \ell(\gamma)\left(
    \lambda+(1-\lambda)\frac{2k\de}r
    \right)
  \le \lambda_0\cdot \ell(\gamma) \le \lambda_0\kappa,
\end{eqnarray*}
as long as $\ell(\gamma)\ge r$. This is true by the
definition of $\gamma$, the assumption on $B(x,r)$ and the
non-contracting derivative of the local diffeomorphisms $T$.
Indeed, if $\ell(\gamma)<r$, then letting $\gamma_1:[0,1]\to
M$ be the only piecewise smooth curve satisfying
$\gamma_1(0)=x$ and $T^n\circ\gamma_1=\gamma$, we must have
$\ell(\gamma)\ge\ell(\gamma_1)$ and thus $\gamma_1(1)\in
B(x,r)$ and $T^n(\gamma_1(1))=y$. This contradicts our
assumption that $y\not\in T^n(B(x,r))$.

Now if we choose $\gamma_0$ such that
$\gamma_0=T^{n-1}\circ\gamma_1$ then, as above, we have both
\[
 \ell(\gamma_0)=\ell(T^{n-1}\circ\gamma_1)\ge r
\qand
\ell(T^{n-1}\circ\gamma_1)\le\lambda_0\cdot\ell(T^n\circ
\gamma_1) = \lambda_0\cdot\ell(\gamma).
\]

Hence, by induction on $k$, we get for every $k=1,\dots,n$ that
\[
r\le \ell(T^{n-k}\circ\gamma_1)\le
\lambda_0^k \cdot\ell(\gamma) \le \lambda_0^k \kappa.
\]

However, this cannot be true for arbitrarily big values of
$n$, since $r,\kappa>0$ are fixed and $\lambda_0\in(0,1)$.
This shows that for every $x\in M$ and $r>0$ there is $n$
such that $T^n(B(x,r))=M$, ending the proof of the lemma.
\cqd

\cpr\label{pr.topmix} Let $\hat T: B\subset \RR^n\to
C^{1+\alpha}(M,M), v\to T+v$, with $B$ a ball around the
origin of $\RR^n$, as defined at the Introduction, where
$T:M\to M$ satisfies the conditions of
Theorem~\ref{th.limitequilibrium}. Then for every $v\in B$,
all $x\in M$ and every given $\ep>0$, there exists $n\in\NN$
such that $T^n_v(B(x,\ep))=M$.  \fpr

\dem
We just have to note that the subset $K$ does not depend on
$v$ for the maps $T_v$ since $DT_v=DT$. Hence if $T$ satisfies the
conditions of Theorem~\ref{th.limitequilibrium}, then every
$T_v$ does also. Thus we can use Lemma~\ref{lem.topmix} with
$T_v$ in the place of $T$.
\cqd



\subsection{Non-degenerate random perturbations}
\label{sec:non-deg-pert}

Here we recall the setting of {\it non-degener\-ate random
  perturbations} as defined in \cite{Ar00}. For a complete
list of propositions and proofs see \cite{Ze03}.

We assume that the family $(\theta_{\ep})_ {\ep > 0}$ of
probability measures on $X$ is such that their supports
have non-empty interior
and
$$
\supp (\theta_\ep) \rightarrow \{t_0\} \quad
\text{when} \quad \ep \rightarrow 0, \quad\mbox{such
  that  } T_{t_0}=T.
$$
For $\omega= (\omega_1, \omega_2, \cdots )\in\Omega$ and
for $n \geq 1$ we set
$$
T^n_{\omega} = T_{\omega_n} \circ \cdots \circ
T_{\omega_1}.$$
Given $x \in M$ and $\omega \in \Omega$ we
call the sequence $(T^n_{\omega}(x))_{n \geq 1}$ a {\it
  random orbit} of $x$.

In what follows we need the map $ \tau_x : X
\rightarrow M, \tau_x(t)= T_t(x)$.

\cde \label{def.nondegpert} We say that $(\hat T,
\theta_\ep)_{\ep > 0}$ is a non-degenerate random
perturbation of $T$ if, for every small enough $\ep$ and
fixed $t_*$ in the interior of $\supp(\th_\ep)$, there
is $\delta_1=\delta_1(\ep)> 0$ such that for all $x \in
M$:
   \begin{enumerate}
   \item $\{T_t(x): t \in \supp(\theta_\ep) \}$ contains
     a ball of radius $\delta_1$ around $T_{t_*}(x);$
   \item $(\tau_x)_{*} \theta_{\ep} $ is
     absolutely continuous with respect to $m.$
   \end{enumerate}
\fde

\cre\label{rmk.noatoms}
We note that $\th_\ep$ cannot have atoms because of the
non-degeneracy condition (2) above.
\fre

We outline some interesting consequences of the
non-degeneracy conditions --- for a proof see~\cite{Ar00}.

\begin{itemize}
\item Any stationary measure $\mu^\ep$ is
  absolutely continuous with respect to $m$.
\item $ \supp (\mu^\ep)$ has non-empty
  interior and $T_t(\supp(\mu^\ep)) \subset
  \supp(\mu^\ep) $ for any $t \in \supp
  (\theta_\ep)$.
 \item $\supp(\mu^\ep) \subseteq B(\mu^\ep)$.
\end{itemize}

Here $B(\mu^\ep)$ is the \emph{ergodic basin} of $\mu^\ep$
\[
B(\mu^\ep)=\left\{x\in M:
\frac1n\sum_{j=0}^{n-1}\varphi(T_\omega^j(x))\to\int\varphi\,
d\mu\mbox{  for all } \varphi\in C(M,\RR) \mbox{ and }
\th^\ep\mbox{-a.e. } \omega\in\Omega \right\}
\]
which by the above properties has positive Lebesgue measure
in $M$.

If $T$ is in the setting of
Theorem~\ref{th.limitequilibrium}, then after
Proposition~\ref{pr.topmix} we deduce that, since the
support of a stationary measure $\mu^\ep$ has non-empty
interior and is forward invariant by $T_{t}$ for any
$t\in\supp\th_\ep$, the support must contain $M$, and so
\emph{there exists only one physical measure} $\mu^\ep$ for
all $\ep>0$, because the support is contained in the basin,
$m\bmod 0$.

These non-degeneracy conditions are not too restrictive
since there always exists a non-degenerate random
perturbation of any differentiable map of a compact manifold
of finite dimension with $X$ the closed ball of radius 1
around the origin of a Euclidean space, see~\cite{Ar00}.
In the setting $M=\TT^n$ with additive noise, as explained
at the Introduction, we have also that $Df_t$ may be
identified with $Df$, so that
$\|(Df_t)^{-1}\|=\|(Df)^{-1}\|, t\in X$, which is very
important in our arguments.


\section{Zero-noise limits are equilibrium measures}
\label{sec:zero-noise-limits}

In what follows we present a proof of
Theorem~\ref{th.limitequilibrium}.  Let $f:M\to M$ be a
local diffeomorphism on a manifold $M$ satisfying the
conditions stated in the above mentioned theorem.  Let also
$\hat f: X\to C^{1+\alpha}(M,M), t\mapsto f_t$ be a
continuous family of maps, where $X$ is a metric space with
$f_{t_0}\equiv f$ for some fixed $t_0\in X$, and
$(\th_\ep)_{\ep>0}$ be a family of probability measures on
$X$ such that $(\hat f,(\th_\ep)_{\ep>0})$ is non-degenerate
random perturbation of $f$.

The strategy is to find a \emph{fixed random generating
  partition} for the system $(\hat f,\th_\ep)$ for \emph{every
small $\ep>0$} and use the absolute continuity of the
stationary measure $\mu^\ep$, together with the
``non-contractive'' conditions on $f$ to obtain (using the
same notations and definitions from
Section~\ref{sec:preliminary-results}) a semicontinuity
property for entropy on zero-noise limits.

\begin{T}
  \label{thm.semicontentropy}
Let us assume that there exists a finite partition $\xi$ of $M$
(Lebesgue modulo zero) which is generating
for random orbits, for every small enough $\ep>0$.

Let $\mu^0$ be a weak$^*$ accumulation point of
$(\mu^\ep)_{\ep>0}$ when $\ep\to0$. If $\mu^{\ep_k}\to\mu^0$
for some $\ep_k\to 0$ when $k\to\infty$, then
\footnote{Cowieson-Young~\cite{CoYo2004} obtain the same
  result for random diffeomorphisms without assuming  the
  existence of a uniform generating partition but need
  either a local entropy condition or that the maps $\hat f$ involved
be of class $C^\infty$.}
\[
\limsup_{k\to\infty} h_{\mu^{\ep_k}}((\hat f, \th_{\ep_k}),\xi)
\le h_{\mu^0}(f,\xi).
\]
\end{T}

The absolute continuity of $\mu^\ep$ and the quasi-expansion
enable us to use a random version of the Entropy  Formula for
endomorphisms (for a more general setting see \cite{LeYo88,BhLi98}).

\begin{T}
  \label{thm.randomPesin}
  If an ergodic stationary measure $\mu^\ep$ for a
  $C^{1+\alpha}$ random perturbation $(\hat f,\th_\ep)$ is
  absolutely continuous and $\int\!\!\int
  \log\|Df_t(x)^{-1}\|\,d\th_\ep(t) d\mu^\ep(x)<0$ for a given
  $\ep>0$, then
\[
h_{\mu^\ep}(\hat f, \th_\ep) =
\int\!\! \int \log|\det Df_t(x)| \, d\mu^\ep(x)
d\th_\ep(t).
\]
\end{T}

Putting Theorems~\ref{thm.semicontentropy}
and~\ref{thm.randomPesin} together shows that
$h_{\mu^0}(f)\ge \int \log|\det Df(x)| \, d\mu^0(x)$,
since $\th_\ep\to\delta_{t_0}$ in the weak$^*$ topology when
$\ep\to0$, by the assumptions on the support of $\th_\ep$ in
subsection~\ref{sec:non-deg-pert}. The reverse inequality
holds in general (Ruelle's inequality~\cite{Ru78}) proving
Theorem~\ref{th.limitequilibrium}.

\subsection{Random Entropy Formula}
\label{sec:rand-pesins-form}

Now we explain the meaning of Theorem~\ref{thm.randomPesin}.

Let $\ep>0$ be fixed in what follows.  The Lyapunov
exponents $ \lim_{n\to\infty}
n^{-1}\log\|Df_\omega^n(x)\cdot v\| $ exist for
$\th^\ep\times\mu^\ep$-almost every $(\omega,x)$ and every
$v\in T_x M\setminus\{0\}$, and are always positive in this
setting.  In fact, since the random perturbations are
additive we have $Df_t=Df$ and, moreover, $\mu^\ep$ is
ergodic, absolutely continuous and $\mu^\ep(K)=0$, so for
$\th^\ep\times\mu^\ep$-almost every $(\omega,x)$
\[
\lim_{n\to\infty}\frac1n\sum_{j=0}^{n-1}
\log \|Df_{\omega(j+1)}(f^j_\omega(x))^{-1}\|
= \int \log\|Df_t(x)^{-1}\| \,d\mu^\ep(x) d\th_\ep(t) < 0.
\]
(Setting $\vfi(t,x)=\log\|Df_t(x)^{-1}\|$ then this is just
the Ergodic Theorem applied to $S:\Omega\times M\to
\Omega\times M$ with $\psi=\vfi\circ\pi$, where
$\pi:\Omega\times M\to X\times M, (\omega,x)\mapsto
(\omega(1),x)$.)  This readily ensures the positivity of the
growth exponent in every direction under random
perturbations because
\[
\log \|(Df_\omega^n(x))^{-1}\|
\le
\sum_{j=0}^{n-1} \log \|(Df_{\omega(j+1)}(f^j_\omega(x)))^{-1}\|.
\]
According to the Multiplicative Ergodic Theorem
  (Oseledets~\cite{Os68}) the sum of the Lyapunov exponents (with
  multiplicities) equals the following limit
  $\th^\ep\times\mu^\ep$-almost everywhere
\[
  \lim_{n\to\infty}\frac1n\log|\det Df_\omega^n(x)| =
  \int\!\!\int\log |\det Df_t(x)|\,d\mu^\ep(x) d\th_\ep(t)>0,
\]
and the identity above follows from the Ergodic Theorem,
since the value of the limit is $S$-invariant, thus constant.

Now {Pesin's Entropy Formula} states that for $C^{1+\alpha}$
maps, $\alpha>0$, with positive Lyapunov exponents
everywhere, as in our setting, the metric entropy with
respect to an invariant measure $\mu^\ep$ satisfies the
relation in Theorem~\ref{thm.randomPesin} if, and only if,
$\mu^\ep$ is absolutely continuous.  In general we
integrate the sum of the positive Lyapunov exponents, see
Liu \cite{Li99} for a proof in the $C^2$ setting. In our
setting of the proof that $\mu^\ep\ll m$ implies the Entropy
Formula is an exercise using the bounded distortion provided
by the H\"older condition on the derivative.


\subsection{Random generating partition}
\label{sec:entr-with-gener}

Here we construct the uniform random generating partition assumed in
the statement of Theorem~\ref{thm.semicontentropy}.
In what follows we fix a weak$^*$ accumulation point $\mu^0$
of $\mu^\ep$ when $\ep\to 0$: there exists $\ep_k\to0$ when
$k\to\infty$ such that $\mu=\lim_k\mu^{\ep_k}$.

To understand how to obtain a generating partition, we need
a preliminary result.

For the following lemma, recall that $K=\{x\in M:
\|Df(x)^{-1}\|=1\}$. In what follows $B(K,\de_0)=\cup_{z\in
  K} B(z,\de_0)$ is the $\de_0$-neighborhood of $K$ and
$\rho>0$ is such that $f\mid B(x,\rho)$ is a diffeomorphism
onto its image and $B(x,\rho)$ is a convex neighborhood for
every $x\in M$, as in
Subsection~~\ref{sec:topological-mixing}.  Using uniform
continuity, we let $\rho_0>0$ be such that for every $x,y\in
M$ and $t\in X$, if $\dist(x,y)<\rho_0$, then
$\dist(f_t(x),f_t(y))<\rho$.

\cle\label{le.expandlength} Let $(f_t)_{t\in X}$ be a family
of maps from a non-degenerate (additive) random perturbation. For any
given $\de_0\in(0,\rho_0)$ there exists $\beta>0$ such that
if $x\in M$ and $y\in M\setminus B(K,\de_0)$ are such that
$\delta_0\le\dist(x,y)\le\rho_0$, then
$\dist(f_t(x),f_t(y))\ge\dist(x,y)+\beta$ for every $t\in
X$.  \fle

\dem Let us assume that $x\in K$, let $y\in M\setminus B(K,\de_0)$ be
such that $\dist(x,y)\in[\delta_0,\rho_0]$ and let
$t\in X$ be fixed. By the choice of $\rho$ there is
a smooth geodesic $\gamma:[0,1]\to M$ with
$\gamma(0)=f_t (x)$ and $\gamma(1)=f_t (y)$ and
$\dist(f_t(x),f_t(y))=\int_0^1 \| \gamma'(s))\| \,
ds<\rho$. In addition, there exists a unique smooth curve
$\gamma_0:[0,1]\to M$ such that $f\circ\gamma_0=\gamma$,
$\gamma_0(0)=x$ and $\gamma_0(1)=y$.

Let us set $b=\|Df_t(y)^{-1}\|=\|Df(y)^{-1}\|<1$ and $K(a)=\{z\in
M:\|Df_t(z)^{-1}\|\ge a\}$ for $a\in(0,1)$. Then there must
be $b_1,b_2\in(b,1)$ with $b_1<b_2$ such that $K(b_1)$ (a
compact set) is in the interior of $K(b_2)$ (recall that
$z\mapsto\|Df_t(z)^{-1}\|$ is continuous and we are assuming
that $x\in K$, that is, $\|(Df_t)^{-1}\|$ assumes the value
1).

We notice that $\|Df_t(z)^{-1}\|<b_1$ for all $z\in K(b_2)\setminus
K(b_1)$ and, moreover, that $\Gamma=\gamma^{-1}(K(b_2)\setminus
K(b_1))$ has nonempty interior, thus positive Lebesgue measure on
$[0,1]$. Then
\begin{eqnarray*}
\dist(f_t(x),f_t(y))
&=&
\int_0^1 \| \dot\gamma(s))\| \, ds = \int_0^1
\|Df_t(\gamma(s))\cdot Df_t(\gamma(s))^{-1} \cdot
\dot\gamma(s))\| \, ds
\\ &\ge&
\frac{1}{b_1}  \int_\Gamma \|Df_t(\gamma(s))^{-1}\cdot
\dot\gamma(s))\| \, ds + \int_{[0,1]\setminus\Gamma}
 \|Df_t(\gamma(s))^{-1}\cdot
 \dot{\gamma}(s))\| \, ds
\\ &>&
\int_0^1 \| \dot\gamma_0(s))\| \, ds
\ge \dist(x,y).
\end{eqnarray*}
If $x\in M\setminus K$, then there exists $b\in(0,1)$ such
that both $x,y\in M\setminus K(b)$ and thus we may take
$\Gamma=[0,1]$ in the calculations above, arriving at the
same sharp inequality.

Hence if $\Delta(\delta_0)=\{(x,y)\in M\times M:
\dist(x,y)\ge\delta_0 \mbox{  and  } y\in M\setminus
B(K,\de_0) \}$, then
\[
h: X \times M\times M\to\RR,\quad
 (t,x,y)\mapsto\dist(f_t(x),f_t(y)) -\dist(x,y)
\]
is positive on the compact set
$X\times\Delta(\delta_0)$, since in the calculations
above $t\in X$ was arbitrary. We just
have to take $\beta=\min h\mid X\times\Delta(\delta_0)$.
\cqd


We now are able to construct a random generating partition under the
conditions of Theorem~\ref{th.limitequilibrium}.

Let us take a finite cover $\{B(x_i,\rho_0/2),
i=1,\dots,\ell\}$ of $M$ by $\rho_0/2$-balls, where
$\rho_0>0$ was already defined.  Since
$\mu^0$ is a probability measure, we may assume that
$\mu^0(\partial \xi)=0$, for otherwise we can replace each
ball by $B(x_i,\gamma\rho_0/2)$, for some $\gamma\in(1,3/2)$
and for all $i=1, \dots, k$.  Now let $\xi$ be the finest
partition of $M$ obtained through all possible intersections
of these balls: $\xi=B(x_1,\gamma\rho_0/2)\vee\dots\vee
B(x_\ell,\gamma\rho_0/2)$. In the following lemma we let
$\rho$ stand for this new radius.

\cre
\label{rmk.0boundary}
The partition $\xi$ is such that all atoms of
$\vee_{j=0}^{n-1} (f_\omega^j)^{-1} \xi$ have boundary
(which is a union of pieces of boundaries of open balls)
with zero Lebesgue measure, for all $n\ge1$ and every
$\omega\in\Omega$. Moreover, since $\mu^0$ is $f$-invariant
and $\mu^0(\partial \xi)=0$, then
$\mu^0(\vee_{j=0}^{n-1} f^{-j} \xi)=0$ also, for all $n\ge1$.  \fre

\cle\label{le.generatingpart}
$
 \bigvee_{k=0}^{+\infty} (f_{\omega}^k)^{-1} \xi = \AA \quad$
  when $ n\to+\infty$
for each $\omega\in\Omega$.
\fle

\dem Let us argue by contradiction assuming that there are
two points $x,y$ such that for some fixed $\omega\in\Omega$:
$\dist(f_\omega^j(x),f_\omega^j(y))\in[\de_0,\rho]$ for some
$\de_0>0$ and $y\in(\bigvee_{i= 0}^{n} f_{\omega}^{-i} \xi)
(x)$ for every $n\ge1$.

Let $\de_1=\min\{\dist(z_1,z_2):z_1,z_2\in K, z_1\neq z_2\}$
be the minimum separation between points in $K$ (we recall that
$K$ is finite) and take
$V=B(K,\min\{\de_1,\de_0\}/4)$. Then it is not possible that
both $f_\omega^j(x),f_\omega^j(y)$ are in the same connected
component of $V$. Using the fact that every
$\rho$-neighborhood is a convex neighborhood and expressing
$\dist(f_\omega^j(x),f_\omega^j(y))$ through the length of a
geodesic, we get a point $z\in M\setminus V$ and
$\be=\be(\de_0,\de_1)>0$ such that
  \begin{eqnarray*}
    \dist(f_\omega^j(x),f_\omega^j(y))
    &=&
    \dist(f_\omega^j(x),f_\omega^j(z))
    +
    \dist(f_\omega^j(z),f_\omega^j(y))
    \\
    &\ge&
    \dist(f_\omega^{j-1}(x),f_\omega^{j-1}(z))
    +
    \dist(f_\omega^{j-1}(z),f_\omega^{j-1}(y))
    + 2\be
    \\ &\geq&
     \dist (f_\omega^{j-1}(x),f_\omega^{j-1}(y)) + 2 \beta
  \end{eqnarray*}
  for every $j>0$, applying Lemma~\ref{le.expandlength}
  twice. But then the upper bound $\rho$ for the distance
  between iterates of $x$ and $y$ cannot hold for all
  $j\ge1$. This shows that $\de_0$ cannot be positive, hence
  the diameter of the atoms of the refined partitions tends
  to zero. This is enough to conclude the statement of the
  lemma. \cqd

This last lemma implies that $\xi$ is a random generating
partition as in the statement of the Random Kolmogorov-Sinai
Theorem~\ref{thm.KSrandom}. Hence we conclude that
$h_{\mu^{\ep_k}}((\hat
f,\theta_{\ep_k}),\xi)=h_{\mu^{\ep_k}}(\hat f,
\theta_{\ep_k})$ for all $k\ge1$.


\subsection{Semicontinuity of entropy on zero-noise}
\label{sec:semic-entr-zero-1}

Now we start the proof of Theorem~\ref{thm.semicontentropy}.
We need to construct a sequence of partitions of
$\Omega\times M$ according to the following result. For a
partition $\cP$ of a given space $Y$ and $y\in Y$ we denote
by $\cP(y)$ the element (atom) of $\cP$ containing $y$.
We set $\omega_0=(t_0,t_0,t_0,\dots)\in\Omega$ in what follows.

\cle\label{lem.seqnpartitions} There exists an increasing
sequence of measurable partitions $(\BB_n)_{n\ge1}$ of
$\Omega$ such that
 \begin{enumerate}
 \item $\omega_0\in\inter \BB_n(\omega_0)$ for all $n\ge1$;
 \item $\BB_n \nearrow \BB$, $\th^{\ep_k} \bmod 0$ for all
   $k\ge1$ when $n\to\infty$;
 \item $\lim_{n \to \infty} H_{\rho} (\xi \mid \BB_n)
 = H_{\rho} (\xi \mid \BB)$ for every measurable finite
 partition $\xi$ and any $S$-invariant probability measure $\rho$.
 \end{enumerate}
\fle

\dem For the first two items we let $\CC_n$ be a finite
$\th_{\ep_k}\bmod 0$ partition of $X$ such that
$t_0\in\inter\CC_n(t_0)$ with $\diam \CC_n\to0$ when
$n\to\infty$. Example: take a cover $(B(t,1/n))_{t\in X}$ of
$X$ by $1/n$-balls and take a subcover $U_1,\dots,U_k$ of
$X\setminus B(t_0,2/n)$ together with $U_0=B(t_0,3/n)$; then
let $\CC_n=U_0\vee\dots\vee U_k$.

We observe that we may assume that the boundary of these
balls has null $\th_{\ep_k}$-measure for all $k\ge1$, since
$(\th_{\ep_k})_{k\ge1}$ is a denumerable family of
non-atomic probability measures on $X$ (see
Remark~\ref{rmk.noatoms}).  Now we set
\[
\BB_n=\CC_n\times\stackrel{n}{\dots}\times\CC_n\times\Omega\quad
\mbox{for all  }n\ge1.
\]
Then since $\diam \CC_n\le 2/n$ for all $n\ge1$ we have that
$\diam \BB_n\le 2/n$ also and so tends to zero when
$n\to\infty$.  Clearly $\BB_n$ is an increasing sequence of
partitions. Hence $\vee_{n\ge1} \BB_n$ generates the
$\sigma$-algebra $\BB$, $\th^{\ep_k} \bmod0$ (see e.g.
\cite[Lemma 3, Chpt. 2]{Bi65}) for all $k\ge1$. This proves
items (1) and (2).

Item (3) of the statement of the lemma is Theorem 12.1 of
Billingsley~\cite{Bi65}.
\cqd

Now we use some properties of conditional entropy to obtain
the right inequalities. We start with
\begin{eqnarray*}
  h_{\mu^{\ep_k}}(\hat f,\theta_{\ep_k})
  &=&
  h_{\mu^{\ep_k}}((\hat f,\theta_{\ep_k}),\xi) = h_{
  \theta^{{\ep_k}} \times \mu^{\ep_k}}^{ \BB \times M} (S,
  \Omega\times\xi)
  \\
  &=&
  \inf\frac1n H_{\th^{\ep_k}\times\mu^{\ep_k}}\left(
  \bigvee_{j=0}^{n-1} (S^j)^{-1}(\Omega\times\xi) \mid
  \BB\times M \right)
\end{eqnarray*}
where the first equality comes from
subsection~\ref{sec:entr-with-gener} and the second one can
be found in Kifer~\cite[Thm. 1.4, Chpt. II]{Ki86}, with
$\Omega\times\xi=\{ \Omega\times A: A\in\xi\}$.
Hence for arbitrary fixed $N\ge1$ and for any $m\ge1$
\begin{eqnarray*}
 h_{\mu^{\ep_k}}(\hat f,\theta_{\ep_k})
 &\le&
 \frac1N H_{\th^{\ep_k}\times\mu^{\ep_k}}\left(
  \bigvee_{j=0}^{N-1} (S^j)^{-1}(\Omega\times\xi) \mid
  \BB\times M \right)
 \\
 &\le&
 \frac1N H_{\th^{\ep_k}\times\mu^{\ep_k}}\left(
  \bigvee_{j=0}^{N-1} (S^j)^{-1}(\Omega\times\xi) \mid
  \BB_m\times M \right)
\end{eqnarray*}
because $\BB_m\times M\subset\BB\times M$. Now we fix $N$ and $m$,
let $k\to \infty$ and note that since
$\mu^0(\partial\xi)=0=\de_{\omega_0}(\partial \BB_m)$ it must be
that
\[
(\de_{\omega_0}\times\mu^0)(\partial(B_i\times\xi_j))=0\quad
\mbox{for all  } B_i\in\BB_m \mbox{  and  } \xi_j\in\xi,
\]
where $\de_{\omega_0}$ is the Dirac mass concentrated at
$\omega_0\in\Omega$.  Thus we get by weak$^*$ convergence of
$\th^{\ep_k}\times\mu^{\ep_k}$ to
$\delta_{\omega_0}\times\mu^0$ when $k\to\infty$
\begin{equation}
  \label{eq:5}
  \limsup_{k\to\infty} h_{\mu^{\ep_k}}(\hat
  f,\theta_{\ep_k})
  \le
  \frac1N H_{\de_{\omega_0}\times\mu^0}\left(
  \bigvee_{j=0}^{N-1} (S^j)^{-1}(\Omega\times\xi) \mid
  \BB_m\times M \right)
  =\frac1N H_{\mu^0} \big(\bigvee_{j=0}^{N-1} f^{-j}\xi \big).
\end{equation}
Here it is easy to see that the middle conditional entropy
of~\eqref{eq:5} (involving only finite partitions) equals
$N^{-1}\sum_i \mu^0(P_i)\log\mu^0(P_i)$, with
$P_i=\xi_{i_0}\cap f^{-1}\xi_{i_1}\cap\dots\cap f^{-(N-1)}
\xi_{i_{N-1}}$ ranging over every possible sequence of
$\xi_{i_0},\dots,\xi_{i_{N-1}}\in \xi$.

Finally, since $N$ was an arbitrary integer,
Theorem~\ref{thm.semicontentropy} follows from~\eqref{eq:5}.
We have completed the proof of the first part of
Theorem~\ref{th.limitequilibrium}.


\section{Existence of a.c.i.m. and stochastic stability}
\label{sec:exist-acim-stoch}

Let $f:M\to M$ be as in the statement of
Theorem~~\ref{th.limitequilibrium}.
The assumptions on $f$ ensure that for every $x\in M$ and
all $v\in T_x M\setminus\{0\}$ we have
\[
\liminf_{n\to\infty}\frac1n\log\|Df^n(x)\cdot v\|\ge0.
\]
Thus the Lyapunov exponents for any given $f$-invariant
probability measure $\mu$ are non-negative. Hence the sum
$\chi(x)$ of the positive Lyapunov exponents of a
$\mu$-generic point $x$ is such that
\begin{equation}
  \label{eq:8}
\chi(x)=\lim_{n\to\infty}\frac1n\log|\det Df^n(x)|
\quad\mbox{and}\quad
\int\chi\, d\mu=\int\log|\det Df|\, d\mu
\end{equation}
by the Multiplicative Ergodic Theorem and the standard
Ergodic Theorem.

Using Theorem~\ref{th.limitequilibrium} we know that there
is only one stationary measure for every $\ep>0$ (see the
beginning of Section~\ref{sec:zero-noise-limits}), and every
weak$^*$ accumulation point $\mu$ of the stationary measures
$(\mu^\ep)_{\ep>0}$, when $\ep\to0$, is an equilibrium state
for $-\log|\det Df|$, that is~~\eqref{eq:3} holds.  We may
and will assume that $\mu$ is ergodic due to the following

\cle\label{lem.ergdecompequil} Almost every ergodic
component of an equilibrium state for $-\log|\det Df|$ is
itself an equilibrium state for this same function.  \fle

\dem Let $\mu$ be an $f$-invariant measure satisfying
$h_\mu(f)=\int \log|\det Df| \, d\mu$. On the one hand,
the Ergodic Decomposition Theorem (see e.g
Mañé~\cite{Man87}) ensures that
\begin{equation}
  \label{eq:9}
\int \log|\det Df| \, d\mu = \int\!\! \int \log|\det
Df| \, d\mu_z \, d\mu(z) \quad\mbox{and}\quad h_\mu(f)=\int
h_{\mu_z}(f) \, d\mu(z).
\end{equation}
On the other hand, Ruelle's inequality guarantees for a
$\mu$-generic $z$ that (recall~\eqref{eq:8})
\begin{equation}
  \label{eq:10}
h_{\mu_z}(f)\le\int \log|\det Df| \, d\mu_z.
\end{equation}
By~\eqref{eq:9} and~\eqref{eq:10}, and because $\mu$ is an
equilibrium state~\eqref{eq:3}, we conclude that we have
equality in~\eqref{eq:10} for $\mu$-almost every $z$.
\cqd

Now we note that since $K$ is finite, if $\mu(K)>0$, then
$\mu$ (which is ergodic) is concentrated on a periodic orbit. Hence
$h_\mu(f)=0$ and so by the entropy formula these orbits are
non-volume-expanding (the Jacobian equals 1).

Finally, for an ergodic equilibrium state $\mu$ with
$\mu(K)=0$, we must have
\[
\lim_{n\to\infty}\frac1n\log\| (Df^n(x))^{-1}\|\le
\lim_{n\to\infty}\frac1n\sum_{j=0}^{n-1}\log\|Df(f^j(x))^{-1}\|
=\int \log\|(Df)^{-1}\| \, d\mu<0 ,
\]
$\mu$-almost everywhere. This means that the Lyapunov
exponents of $\mu$ are strictly positive, where $f$ is a
$C^{1+\alpha}$ endomorphism, $\alpha>0$.

Now the extension of the Entropy Formula for endomorphisms
obtained by Liu~\cite{Li98} (in the $C^2$ setting, but the
distortion estimates need only a H\"older condition on the
derivative) ensures that an equilibrium state whose Lyapunov
exponents are all positive must be absolutely continuous
with respect to Lebesgue measure.  Hence $\mu\ll m$.

The previous discussion shows that $f$ is non-uniformly
expanding in the sense of Alves-Bonatti-Viana~\cite{ABV00}.
We obtain that
\begin{equation}
  \limsup_{n\to+\infty}\frac1n\sum_{j=0}^{n-1}
  \log\|Df(f^j(x))^{-1}\| < 0, \quad \mu-\mbox{almost every  }
  x.
\end{equation}
These authors show that any absolutely continuous invariant
measure $\mu$ in this setting has basin $B(\mu)$ containing
an open subset $U$ Lebesgue modulo 0. By the topological
mixing property of $f$ (Lemma~\ref{lem.topmix}) and because
$f$ is a regular map, we deduce that $B(\mu)$ must contain
all of $M$, Lebesgue modulo 0, and is thus unique.

These arguments hold true for every ergodic component of any
weak$^*$ accumulation point of stationary measures when
$\ep\to0$, thus every such accumulation point is a linear
convex combination of an absolutely continuous invariant
measure with finitely many Dirac masses concentrated on
non-volume-expanding orbits.  This ends the proof of
Theorem~\ref{th.limitequilibrium}.

\subsection{Stochastic stability}
\label{sec:stability}

Now we assume that $f:M\to M$ satisfies all the conditions
of the statement of Theorem~\ref{th.limitSRB}. We arrive at
the same conclusion of Theorem~\ref{th.limitequilibrium} but
since we assume that $|\det Df|>1$ on $K$, we would arrive
at a contradiction if $\mu$ is an equilibrium state with
$\mu(K)=0$, that is, the expanding volume condition avoids
Dirac masses on periodic orbits as ergodic components of
zero-noise limit measures. Thus in this setting we must have
$\mu(K)=0$ for all ergodic equilibrium states. Hence there
is only one equilibrium state: the unique absolutely
continuous invariant probability measure $\mu$ of $f$.

Therefore every weak$^*$ accumulation point of the stationary
measures, when $\ep\to0$, must be equal to $\mu$, showing the
stochastic stability of $\mu$ and concluding the proof of
Theorem~\ref{th.limitSRB}.


\section{Random perturbations of the intermittent map}
\label{sec:rand-pert-interm}

Let $T$ be defined
as in Introduction and consider a non-degenerate random
perturbation of $T$.  For any small $\ep>0$, we know that
there exists a single stationary measure $\mu^\ep$, see
Subsection~\ref{sec:topological-mixing}.

\subsection{Characterization of zero-noise measures }

Let $T$ be as above for $\alpha \in (0,1)$.
Here we prove Theorem \ref{th.equil} in two steps.
Let $\mu$ be a weak$^*$ accumulation point of
$\mu^{\ep}$, when $\ep\to0$. We show that $\mu \in \EE = \{ t
\delta_0 + (1-t) \mu_{\rm SRB} : 0 \leq t \leq 1 \}$.

Firstly we prove that $ h_{\mu}(T) = \int \log DT \, d \mu$,
and in the sequel we deduce that any $T$-invariant measure
satisfying the entropy formula as above should belong to
$\EE.$ The latter is proved also in \cite{Pol} by different methods.

As we are considering additive random perturbation of $T$,
we can apply Theorem~\ref{th.limitequilibrium} and
conclude that $\mu$ is an equilibrium state of $-\log
DT.$ As in the previous section, we may and will assume
that $\mu$ is ergodic, by Lemma~\ref{lem.ergdecompequil}.

Now we consider two cases: either $\mu (\{0\}) > 0$, or $\mu
(\{0\}) = 0$.

In the first case, the ergodicity of $\mu$ ensures that  $\mu =
\delta_0$, since $0$ is a fixed point for $T$.

For the second case, as $\log DT > 0$ for all points in
$\SS$ except $0$, we have that
\[
h_{\mu}(T) = \int \log DT \,d\mu > 0
\]
and the Ergodic Theorem together with the fact that $\SS$ is
one-dimensional guarantees that the Lyapunov exponent of
$\mu$ is positive. Thus $\mu$ is an ergodic probability
measure with positive Lyapunov exponent and positive entropy
which satisfies the Entropy Formula. Now we apply a version
Pesin's Entropy Formula obtained by
Ledrappier~\cite{le81}, which holds for $C^{1+\alpha}$
endomorphisms of $\SS^1$, to conclude that $\mu$ must be
absolutely continuous with respect to Lebesgue measure.

The above arguments show that any typical ergodic
component of every zero-noise limit measure $\mu$ equals either
$\delta_0$ or $\mu_{SRB}$. Hence a straightforward
application the Ergodic Decomposition Theorem to $\mu$
concludes the proof of Theorem \ref{th.equil}.


\subsection{Stochastic stability without a.c.i.m.}
\label{sec:stochastic-stability}

After Theorem \ref{th.limitequilibrium}, any zero-noise
limit measure $\mu$ for the additive random perturbation of
the intermittent map $T$ is an equilibrium state for $-\log
DT$. For $\alpha \geq 1$ this is enough to deduce stochastic
stability of $T=T_\alpha$.

Indeed, let $\mu$ be a weak$^*$ accumulation point of
$\mu^{\ep}$ when $\ep\to0$.  As in the proof of Theorem
\ref{th.equil} (in the previous subsection), we consider the
ergodic decomposition of $\mu$.

We claim that almost all ergodic components of $\mu$ equal
the Dirac measure $\de_0$ concentrated on $0.$ Arguing by
contradiction, we suppose that for some ergodic component
$\eta$ of $\mu$ we have $\eta (\{0\})=0$. Thus by the same
reasoning of the previous subsection (using positive
Lyapunov exponents and Entropy Formula for one-dimensional
maps), this implies that $\eta$ is an absolutely continuous
invariant probability measure for $T.$

However, because for $\alpha \geq 1$ the intermittent map $T$ is
$C^2$, it is well known that $T$ does not admit any absolutely
continuous invariant probability measure in this setting --- see
e.g.~\cite\cite{Vi97b} for a proof of this fact.

Hence if some ergodic component $\eta$ of $\mu$ is such that
$\eta (\{0\})=0$ we arrive at contradiction. Thus
$\eta(\{0\})>0$ and $\eta=\de_0$ by ergodicity, for every
ergodic component of $\mu$. Therefore $\mu=\de_0$.

This proves \emph{stochastic stability of the intermittent
  map when it does not admit absolutely continuous invariant
  probability measures} and ends the proof of
Theorem~\ref{th.stocstable}.

\subsection{Stochastically unstable random perturbation}
\label{sec:unstable}

Here we prove Theorem~\ref{thm.non-stochstable}.
Let $f_t$ be the family \eqref{eq:nonstochfamily} and
let $0<s<1$. Then there exists a unique fixed source
\[
p_s=\frac12\left( \frac{s-1}{s-2} \right)^{1/\alpha}
\in(0,1/2)
\quad\mbox{such that}\quad
f_s'(p_s)=1+\alpha(1-s)>1.
\]
Now we choose $u\in(s,1)$ such that $f_t'\mid [p_u,p_s] >1$
for all $t\in[s,u]$. For this we just have to take $u$ close
enough to $s$.

Clearly $f_u^n(x)\to0$ for all $x\in(1-p_u,p_u)$ when
$n\to\infty$, see Figure~\ref{fig1} and recall that the maps
$f_t$ are symmetric ($f_t(x)=1-f_t(1-x)$) on
$\SS=[0,1]/0\sim 1$.

\cle\label{le.moveto0} For every $x\in(1-p_u,p_u)$
and every sequence $\un t\in[s,u]^\NN$ we have that $f_{\un
  t}^n(x)\to0$ when $n\to\infty$.  \fle

\dem It is straightforward to check that the graph of
$f_t\mid[0,1/2]$ is below the graph of $f_u$ and above the
graph of $f_s$ for every $t\in(s,u)$. Hence $f_{\un
  t}^n(x)\le f_u^n(x)\to0$ when $n\to\infty$ for every
$x\in(0,p_u)$. Using the symmetry we arrive at $f_{\un
  t}^n(x)\to0$ when $n\to\infty$ for all $x\in(1-p_u,0)$.
\cqd

Now we let $\th$ be any probability measure with support contained in
$[s,u]$ and set $\th_0=\th^\NN$.

\cpr\label{pr.getout} For $\th_0\times m$-almost every $(\un
t,x)\in[s,u]^\NN\times [p_u,1-p_u]$ there exists $n\ge1$ such
that $f^n_{\un t}(x)\not\in[p_u,1-p_u]$.  \fpr

Combining the two results above we conclude that for
$\th_0\times m$-almost every $(\un t,x)\in[s,u]^\NN\times
\SS$ we have that
\[
\frac1n \sum_{j=0}^{n-1}\de_{f^j_{\un t}(x)}\to\de_0
\mbox{  in the weak$^*$ topology when  }
n\to\infty,
\]
finishing the proof of Theorem~\ref{thm.non-stochstable}.

To prove Proposition~\ref{pr.getout} we need the following
result whose proof follows standard steps, using the uniform
expansion and the $C^{1+\alpha}$ condition on every $f_t$.
Let us fix $\un t\in[s,u]^\NN$ and a point $0<r<p_u$ such
that $f_t'(s)>1$ for all $t\in[s,u]$. Let also
$\beta_1=\min\{f_t'(x): x\in[r,1-r], t\in[s,u]\}>1$
and
$\beta_2=\max\{f_t'(x):x\in\SS, t\in[s,u]\}>1$.

\cle\label{le.boundist} There exists $C>1$ such that for any
interval $I\subset[r,1-r]$, all $\un t\in[s,u]^\NN$ and
$k\ge1$ such that $f^j_{\un t}(I)\subset[r,1-r]$ for every
$j=0,\dots,k-1$, it holds
\[
\frac1C\le\frac{(f^k_{\un t})^\prime(x)}{(f^k_{\un
    t})^\prime(y)}\le C
\]
for all $x,y\in I$.
\fle

\begin{proof}[Proof of the Proposition]
For $\un t\in[s,u]^\NN$ we define $E_k(\un t)=(f_{\un
  t}^k)^{-1}[p_u,1-p_u]$ for $k\ge1$. We will show that
$\cap_{k\ge1} E_k(\un t)$ has zero Lebesgue measure for any
$\un t$, which is enough to conclude the statement of the
lemma. In fact, this means that $n(\un t,
x)=\min\{k\ge1:f^k_{\un t}(x)\not\in[p_u,1-p_u]\}$ is finite
for every $x$ in a set $X(\un t)$ with $m(X(\un t))=1$, for
every given $\un t\in[s,u]^\NN$. Thus $\Delta=\cup_{\un
  t\in[s,u]^\NN} \{\un t\}\times X(\un t)$ is measurable and
$(\theta_0\times m) (\Delta)=1$.

Let us fix $\un t\in[s,u]^\NN$, take a nonempty interval $I\subset[p_u,1-p_u]$
and show that $I\cap \cap_{k\ge1} E_k(\un t)$  has zero Lebesgue measure.

Let $k>1$ be the first time such that $f^j_{\un
  t}(I)\not\subset[r,1-r]$. There exists such $k$ since by uniform
  expansion $m(f^k_{\un t}(I))\ge\beta_1^k m(I)$ whenever $f^j_{\un
  t}(I)\subset[r,1-r]$ for $j=0,\dots,k-1$. Now there are two
  possibilities: either $f^k_{\un t}(I)\subset(1-p_u,p_u)$ or we have
  $f^k_{\un t}(I)\cap[p_u,1-p_u]\neq\emptyset\neq[1-r,r]\cap f^k_{\un
  t}(I)$.

In the former case we conclude that $I\cap \cap_{k\ge1}
E_k(\un t)=\emptyset$, and the argument ends. 

In the latter case, we let $F=f^k_{\un t}(I)\cap(1-p_u,p_u)$
and observe that either $F\supset[s,p_u]$ or
$F\supset[1-p_u,1-s]$, so $m(F)\ge p_u-s$. Since $f_{\un
  t}^{k-1}(I)\subset[r,1-r]$ we have $m(f^k_{\un
  t}(I))\le\beta_2(1-2r)$ and hence
\[
\frac{m(G)}{m(I)}\ge C\frac{m(F)}{m(f^k_{\un t}(I))} \ge
C\frac{p_u-r}{\beta_2(1-2r)}
\]
where $G=(f^k_{\un t}\mid I)^{-1}({\rm closure\,}(F))$ and
$C>0$ is a bounded distortion constant from Lemma~\ref{le.boundist}.

This shows that $m(I\setminus G)\le(1-C(p_u-r)/(\beta_2(1-2r)))m(I)$.
We may take $r$ so close to $p_u$ that
\[
0<\gamma=1-C\frac{p_u-r}{\beta_2(1-2r)} <1.
\]
If $m(I\setminus G)=0$, we are done. Otherwise we apply the same
argument to each connected component of $I_1=I\setminus G$
inductively, as follows.

Let $I_k\subset I$ be a compact set formed by finitely many
pairwise disjoint closed intervals $I_k=I_{k,1}\cup\dots\cup
I_{k,i_k}$ such that for each $j\in\{1,\dots,i_k\}$ there is
a maximal iterate $n_j$ so that $f_{\un
  t}^n(I_{k,j})\subset[r,1-r]$ for all $n=1,\dots,n_j-1$.

We observe that $\gamma$ does not depend on the number of
iterates of the first exit. Then either $f^{n_j}_{\un
  t}(I_{k,j})\subset(1-p_u,p_u)$, or there exists
$G_{k,j}\subset I_{k,j}$ maximal such that $f^{n_j}_{\un
  t}(G_{k,j})\subset(1-p_u,p_u)$ and $m(I_{k,j}\setminus
G_{k,j})\le\gamma\cdot m(I_{k,j})$, as before.

In the former case we delete $I_{k,j}$ from $I_{k+1}$. In
the latter case, we add the connected components $I^\prime$
of $I_{k,j}\setminus G_{k,j}$ to $I_{k+1}$ and associate to
each of them the maximal number of iterates $n^\prime$ such
that $f^j_{\un t}(I^\prime)\subset[r,1-r]$ for
$j=1,\dots,n^\prime-1$.  Then $m(I_{k+1})\le\gamma \cdot m(I_k)$.
This shows that $m(I_k)\le\gamma^k m(I)\to0$ when
$k\to\infty$. Since by construction $\cap_{k\ge1} I_k$ contains
$I\cap \cap_{k\ge1} E_k(\un t)$, the proof is complete.
\cqd

\end{document}